\documentclass[12pt]{article}

\usepackage{amsmath}
\usepackage{theorem}
\usepackage{amssymb}
\usepackage{latexsym}

\theoremstyle{plain}
\newtheorem{thm}{Theorem}[section]
\newtheorem{lem}[thm]{Lemma}
\newtheorem{prop}[thm]{Proposition}

\newtheorem{df}[thm]{Definition}
\newtheorem{rem}[thm]{Remark}
\newenvironment{pf}{\noindent{\em Proof.}}{\hbox{}\hfill $\Box$\vskip0.2cm}
{\theorembodyfont{\rmfamily} \newtheorem{exa}[thm]{Example}}

\newcommand{\Q}{\mathbb{Q}}

\newcommand{\Z}{\mathbb{Z}}

\newcommand{\mf}[1]{\mathfrak{#1}}

\newcommand{\ad}{\mathop{\mathrm{ad}}}
\newcommand{\GL}{\mathop{\mathrm{GL}}}

\newcommand{\Sym}{\mathop{\mathrm{Sym}}}
\newcommand{\Der}{\mathop{\mathrm{Der}}}
\newcommand{\Stab}{\mathop{\mathrm{Stab}}}
\newcommand{\gl}{\mathfrak{\mathop{gl}}}
\newcommand{\ssl}{\mathfrak{\mathop{sl}}}

\newcommand{\Ad}{\mathop{\mathrm{Ad}}}

\newcommand{\PGL}{\mathop{\mathrm{PGL}}}
\newcommand{\Aut}{\mathop{\mathrm{Aut}}}
\newcommand{\Pic}{\mathop{\mathrm{Pic}}}
\newcommand{\Div}{\mathop{\mathrm{Div}}}

\newcommand{\p}{\mathbb{P}}

\newcommand{\PxP}{\mathop{\mathbb{P}^1\times\mathbb{P}^1}}
\newcommand{\slxsl}{\mathop{\mathfrak{sl}_2\oplus\mathfrak{sl}_2}}

\title{Parametrizing Del Pezzo surfaces of degree 8 using Lie algebras}

\author{Willem A. de Graaf, University of Trento, Italy \and 
        Jana P\'\i lnikov\'a, RICAM, Linz, Austria\protect \and
        Josef Schicho, RICAM, Linz, Austria}

\begin{document}

\maketitle

\begin{abstract}
  For a Del Pezzo surface of degree 8 given over the rationals we
  decide whether there is a rational parametrization of the surface
  and construct one in the affirmative case. We define and use 
  the Lie algebra of the surface to reach the aim.
  The algorithm has been implemented in Magma.
\end{abstract}

\section{Introduction}

In this paper we decide whether a given Del Pezzo surface of degree 8
has a rational parametrization over $\Q$, and find one in the affirmative
case. There are two kinds of Del Pezzo surfaces of degree 8: blowups
of $\mathbb{P}^2$ in one point and twists of $\mathbb{P}^1\times\mathbb{P}^1$.
Here we deal with both of them.

The problem is a particular case of a more general one, namely parametrization
of surfaces over the rational numbers. There one reduces to several base
cases. Except some trivial cases (e.g. $\mathbb{P}^2$), there are Del Pezzo
surfaces of degrees 5 till 9, and conic fibrations. The latter are solved
in \cite{josef1}. Rational parametrization of Del Pezzo surfaces
of degree~5 is discussed in \cite{Shepperd-Barron:92}, and of degree~7
can be found in e.g. \cite{manin}. Parametrization of degree~9 is
solved in~\cite{ghps}. Hence the last unsolved class are Del Pezzo surfaces
of degree~6.

For Del Pezzo surfaces of degree~6 the Hasse principle holds
(cf \cite{manin}). For a given prime $p$ we can decide (e.g. using
Magma~\cite{magma}) whether there is a point on the surface over the local
field $\Q_p$. But there still remains the problem of finding a finite
set of ``bad primes'' and finding a rational point provided that we have
points over local fields.

Here we use another approach. Similarly as in \cite{ghps}, the 
parametrization problem is reduced to a problem concerning Lie algebras
and their representations. 
The Lie algebra approach appears to be preferable even in situations where
also other methods are known, such as the parametrization problem
of blowups of the plane in one point, see~\cite{manin}.

The paper is structured as follows. In section~\ref{sec:la},
we introduce the Lie algebra of a variety and give an algorithm
for computing it. Section~\ref{sec:id} is dealing with ``identification
problems''; the main focus will be to reduce the problem of
identifying a variety (i.e. constructing an isomorphism if exists)
to the problem of identifying a Lie algebra. This will solve
some instances of our parametrization problem of Del Pezzo surfaces
of degree~8. The remaining instances are dealt with in 
section~\ref{sec:simple}.

The second and third author were supported by SFB Grant F1303
of the Austrian Science Fund (FWF).

\section{The Lie Algebra of a Variety} \label{sec:la}

In this section, we define the Lie algebra of a variety and give
a method for computing it.

Throughout, we assume that $F$ is a field of characteristic zero. We are mostly
interested in the case $F=\Q$, but the algorithm to be described 
works equally well for number fields, and there is one step where
field extensions are needed (see section~\ref{sec:simple}).

Let $X$ be a projective variety over $F$.
We denote the group of its automorphisms by $\Aut(X)$. 
The first idea to define the Lie algebra of $X$ would be to take
the tangent space of $\Aut(X)$ at the identity,
but this does not work in general because $\Aut(X)$ need not be
an algebraic group. Hence we introduce
\[  {\Aut}_0(X) = \left\{\varphi\in\Aut(X) \mid
  \varphi \textrm{ acts trivially on }\Pic(X)\right\}. \]
The advantage of working with $\Aut_0(X)$ rather than
with $\Aut(X)$ become clear by this

\begin{thm}
The group $\Aut_0(X)$ is an algebraic group over $F$. 
For any very ample divisor $D$ of $X$, there is a faithful representation
of $\Aut_0(X)$ into $\PGL_{n+1}(F)$, where $n:=\dim(D)$.
\end{thm}

\begin{pf}
Let $i:X\hookrightarrow\mathbb{P}^n$ be the
rational map associated to $D$. It is an embedding since $D$ is very ample.
Suppose that $\varphi\in\Aut_0(X)$. Then the pullback
$\varphi^*:\Div(X)\to\Div(X)$ transforms the complete linear system
$|D|$ into itself. Routine calculation shows that this transformation
is projective and its dual is an extension of $\varphi$ to the ambient
projective space $\p^n$. Clearly, $X$ contains $n+2$ points in
general position (not necessarily defined over $F$), hence any 
projective transformation leaving $X$ pointwise fixed is the identity,
and it follows that the representation of $\Aut_0(X)$ into $\PGL_{n+1}(F)$
is faithful.

If we choose a different very ample divisor $D_1$, then the composition
of homomorphisms to and from $\Aut_0(X)$ is algebraic (i.e. it
is a regular function). Hence the the algebraic structure of $\Aut_0(X)$ 
does not depend on the choice of $D$.
\end{pf}

For convenience, we would prefer a representation in $\GL_{n+1}(F)$
rather than in $\PGL_{n+1}(F)$. This is not always possible, for instance,
it is not clear how to embed $\Aut_0(\p^1)$ into $\GL_2(\Q)$.
But after passing to Lie algebras the situation is much easier:
we just add a direct one-dimensional abelian summand in order to compensate
the difference between $\gl_{n+1}(F)$ and $\ssl_{n+1}(F)$, the Lie algebra of
$\PGL_{n+1}(F)$.

\begin{df}
Let $X$ be a projective variety over $F$.
We define $L_0(X,F)$ as
the Lie algebra of the algebraic group $\Aut_0(X)$, and
$L(X,F)$ as the direct sum of $L_0(X,F)$ and the abelian one-dimensional
Lie algebra $C$.

We also write $L_0(X)$ and $L(X)$ as shorthands, if there
is no ambiguity of the field. We refer to $L_0(X)$ as
{\em the Lie algebra of the variety $X$}.
\end{df}

Here is a theorem that can be used for the computation of $L(X,F)$.

\begin{thm}
Let $X$ be a projective variety such that $\Pic(X)$ is discrete.
Let $D\in\Div(X)$ be a very ample divisor, and let $n:=\dim(D)+1$. 
Let $i:X\to\p^{n-1}$ be the associated embedding.
Let $\Aut_i(X)\subset\GL_n(F)$ be the group of all invertible linear maps
whose projectivization maps $i(X)$ into itself. Then
$L(X,F)$ is the Lie algebra of $\Aut_i(X)$.
\end{thm}

\begin{pf}
Note that $\Aut_i(X)$ is an algebraic group, because it can be given
by polynomial equations, namely $g\in\Aut_i(X)$ if and only if $f_i(gp) = 0$ 
for all $p\in X$ and all $i$ such that $f_i$'s generate the vanishing ideal
of the embedded variety $i(X)$.
The multiplicative group $Z$ of scalar matrices is an algebraic subgroup
in the center of $\Aut_i(X)$.
The quotient group $\Aut_i(X)/Z$ is an algebraic group of automorphisms of $X$
containing $\Aut_0(X)$. Because the Picard group is discrete,
the connected component of the identity of $\Aut_i(X)/Z$ leaves it
pointwise fixed, hence it is contained in $\Aut_0(X)$.
It follows that $\Aut_i(X)/Z$ and $\Aut_0(X)$ have the same
Lie algebra, namely $L_0(X,F)$.

Because $Z$ is contained in the center of $\Aut_i(X)$, its Lie algebra $C$
is contained in the center of the Lie algebra of $\Aut_i(X)$. 
It follows that $C$ is a direct summand,
and the co-summand is the Lie algebra of the quotient.
\end{pf}

Note that by this construction, $L(X,F)$ is a Lie algebra of matrices,
so the construction gives not only $L(X,F)$ but also a representation
$L(X,F)\hookrightarrow\gl_n(F)$; of course, the representation 
depends on the embedding $i$.

\begin{exa} \label{exa:p}
Let $r>0$. Let $X = \p^r$. 
Then every automorphism fixes the Picard group, which is isomorphic to $\Z$. 
So we have $\Aut(X) = \Aut_0(X) = \PGL_{r+1}(F)$,
and $L_0(X,F)=\ssl_{r+1}(F)$, and $L(X,F)=\gl_{r+1}(F)$.

Let $d>0$. Let $D$ be a divisor of degree $d$. Then
$n=\dim(D)={r+d+1\choose d}-1$, and the associated map $i:X\to\p^n$
is the $d$-uple embedding. The group $\Aut_i(X)$ is the
$d$-th symmetric power of $\GL_{r+1}(F)$, and its Lie algebra
is the representation of $\gl_{r+1}(F)$ by $d$-th symmetric powers.
\end{exa}

The paper \cite{ghps} contains the following converse of 
Example~\ref{exa:p}: if $X$ is a twist of $\p^r$ and
$L(X,F)\cong \gl_{r+1}(F)$, then $X\cong \p^r$. The
problem of constructing the isomorphism from $\p^r$ to
$X$ can be reduced to the construction of a Lie algebra
isomorphism from $\gl_{r+1}(F)$ to $L(X,F)$. In section~\ref{ss:pxp},
we will prove a similar result for twists of $\p^1\times\p^1$.

In the applications we are interested, the ideal of the variety $X$ can be
given by quadratic equations. This is equivalent to $D$ having
the property $N_1$ (see \cite{Schenck:04}).
In this case, there is a particularly easy way of computing its Lie algebra.

\begin{thm}\label{thm:alglie}
Let $X\subset\p^n$ be an embedded projective variety.
Assume that the ideal of $X$ is generated by quadrics.
Write all of these quadratic equations as $p^TAp$,
where $A$ is a symmetric matrix of size $(n+1)\times(n+1)$.
Let $I$ be the linear space generated by these matrices.

Then the Lie algebra $L(X,F)$ is the matrix algebra
\[ \{ x\in \gl_{n+1}(F) \mid x^TA+Ax\in I \text{ for all }A\in I\} . \]
\end{thm}

\begin{pf}
Let $i$ be the embedding of $X$. We have
\[ {\Aut}_i(X) = \{ g\in {\GL}_{n+1}(F) \mid g^TAg\in I\text{ for all }A\in I\} .\]
Let $W$ denote the vector space of $n+1\times n+1$-matrices over $F$. We 
have a rational representation $\rho : \GL_{n+1}(F)\to \GL(W)$ given by
$\rho(g)(A) = g^TAg$. Then ${\Aut}_i(X)$ is the group of all $g\in \GL_{n+1}(F)$
such that $\rho(g)I=I$. By \cite{chevalley}, Corollary 1 to Theorem 1,
Chapter III, No 9, the Lie algebra of ${\Aut}_i(X)$ consists of all 
$x\in \gl_{n+1}(F)$ such that $(d\rho)(x)(I)\subset I$. Now $(d\rho)(x)(A)
= x^TA+Ax$. 
\end{pf}

Of course, it is sufficient to collect all conditions for $A$ in a fixed
basis of $I$. Hence $L(X,F)$ can be computed by linear algebra.

\section{Identification Problems} \label{sec:id}

In this section we treat some special instances and subproblems of
the pa\-ra\-me\-tri\-za\-tion problem for Del Pezzo surfaces of degree~8,
of the following type: given a variety $X$, decide whether it is equivalent
to a fixed variety $Y$; and if yes, construct an isomorphism from
$Y$ to $X$. We call this type of problem the {\em identification problem}
for $Y$.

We solve it for $Y=\p^1,\PxP$, and the blowup of $\p^2$ at a single point.
We denote this blowup variety by $\mathbb{Y}$.
(For $Y=\p^2$, the problem was solved in \cite{ghps}.)
In all these cases, the existence part of the identification problem
is not difficult to solve when $F$ is algebraically closed. So we
may assume that $X$ is a twist of $Y$, i.e. that $X$ is isomorphic
to $Y$ over the algebraic closure of $F$.

In all the cases above, the anticanonical divisor $-K$ is very ample.
If $X$ and $Y$ are isomorphic, then the anticanonical images
$i_{-K}(X)$ and $i_{-K}(Y)$ are projectively isomorphic, because
any isomorphism induces a linear isomorphism between the spaces of
global sections of the two anticanonical line bundles (see also \cite{ghps}).
The problem of deciding whether two embedded projective varieties
are projectively isomorphic, and to construct a projective
transformation if exists, will be called the {\em embedded identification
problem}.

A necessary condition for $X$ being isomorphic to $Y$ is that
$L_0(X)$ is isomorphic to $L_0(Y)$. If both $X$ and $Y$ are
anticanonically embedded, the isomorphism $Y\to X$ is described
by $p\mapsto Mp$ for some matrix $M\in\GL_{n+1}(F)$, 
where $n=\dim(-K_X)=\dim(-K_Y)$.
Then we also have a Lie algebra isomorphism
$\nu:L_0(Y)\to L_0(X)$ given by $\nu(x)=MxM^{-1}$ for the same matrix $M$.
The matrix $M$ defines an isomorphism of the $L_0(Y)$-modules
given by the inclusion $L_0(Y)\hookrightarrow \gl_{n+1}(F)$ and
by its composition with the Lie algebra isomorphism $\nu$.

Given $X$, we claim that the embedded identification problem for $Y$ 
can be solved by the following algorithm (assuming $Y$ is one of
$\p^1$, $\PxP$, or $\mathbb{Y}$):
\begin{enumerate}
\item Compute $L_0(X)$ and solve the {\em Lie algebra identification problem}
        for $L_0(Y)$; i.e., construct a Lie algebra isomorphism $\nu$
        if exists. Otherwise, $X$ and $Y$ are not isomorphic.
\item Construct an isomorphism $M$ between the $L_0(Y)$-modules
        defined by the inclusion $L_0(Y)\hookrightarrow \gl_{n+1}(F)$
        and by the composition with $\nu$. If the modules are
        not isomorphic, then $X$ and $Y$ are not projectively
        equivalent.
\item Check if $M$ transforms $Y$ to $X$. If yes, we have found the
        isomorphism. Otherwise, $X$ and $Y$ are not projectively
        equivalent.
\end{enumerate}

Methods for solving the Lie algebra identification problem (step~1)
and for solving the module identification problem (step~2)
will be explained in the subsequent subsections.

The correctness of the algorithm follows from the following statements.
\begin{itemize}
\item Assume $Y$ and $X$ are projectively equivalent via some matrix $M$, 
        and $\nu:L_0(Y)\to L_0(X)$ is a Lie algebra isomorphism.
        Then conjugation by $M$ is another isomorphism from $L_0(X)$ 
        to $L_0(Y)$ (cf. \cite{ghps}, Proposition 3.4).
        Composing these two, we get a Lie algebra
        automorphism of $L_0(Y)$. By Lemma~\ref{lem:lift} below,
        this automorphism is equal to the conjugation by a matrix 
        $N\in \Aut_i(Y)$.
        Then $\nu$ is equal to conjugation by $NM$, and
        $NM$ is an isomorphism of the two $L_0(Y)$-modules
        in step~2. In particular, these $L_0(Y)$-modules are 
        isomorphic.
\item Assume $Y$ and $X$ are projectively equivalent via some matrix $U$, and 
        conjugation by $M\in\GL_{n+1}(F)$ is a homomorphism from
        $L_0(X)$ into $L_0(Y)$
        (this is equivalent to $M$ being a module isomorphism as
        computed in step~2). Then we claim that $M$ maps $Y$ to $X$.
        For this we may assume that the ground field is algebraically
        closed. Indeed, by Theorem \ref{thm:alglie}, $L_0(Y,\overline{F})=
        L_0(Y,F)\otimes \overline{F}$, and similarly for $L_0(X,\overline{F})$
        Furthermore, conjugation by $M$ is also a homomorphism of 
        $L_0(X,\overline{F})$ into $L_0(Y,\overline{F})$. Note that
        conjugation by $U^{-1}M$ is an
        automorphism of $L_0(Y)$. By Lemma~\ref{lem:lift} below,
        $U^{-1}M$ transforms $Y$ to itself, hence $M=U(U^{-1}M)$
        transforms $Y$ to $X$.
\end{itemize}

\begin{lem}\label{lem:surjective}
Let $G\subset \GL(V)$, $H\subset \GL(W)$ be algebraic groups, where
$V$ and $W$ are vector spaces over an algebraically closed field.
Assume that $G$ and $H$ have the same dimension, and the same number of
connected components. Let $\sigma : G\to H$ be an injective rational
representation of $G$. Then $\sigma$ is surjective.
\end{lem}

\begin{pf}
Let $G^0,H^0$ be the connected components of the identity. Then
$\sigma(G^0)$ is connected as well. Hence $\sigma(G^0)\subset H^0$.
Furthermore, by \cite{chevalley}, Corollary 1 to Proposition 2,
Chapter II, \S 7, $\sigma(G^0)$ is an algebraic subgroup of $H$.
Furthermore, $\dim \sigma(G^0)=\dim G^0$ (this follows for example
from \cite{chevalley}, Corollary to Proposition 8, Chapter II, 
\S 6). Hence $\sigma(G^0)$ has finite index in $H^0$ (\cite{chevalley},
Proposition 4, Chapter II, \S 6). Since $H^0$ is the unique irreducible 
algebraic subgroup of $H$ of finite index, we conclude that $\sigma(G^0)=H^0$.
Let $x_1G^0,\ldots,x_tG^0$ be the irreducible components of $G$. Then
$\sigma(x_iG^0)= \sigma(x_i)H^0$. Hence the $\sigma(x_i)H^0$ are
irreducible components of $H$. Since $H$ has the same number of such 
components as $G$, we conclude that $\sigma$ is surjective.
\end{pf}

\begin{lem} \label{lem:lift}
Suppose that the ground field is algebraically closed, and 
that the centre of $L_0(Y)$ is $0$. 
Let $Z=\{\lambda I_{n+1}\}$, where $I_{n+1}\in \GL_{n+1}(F)$ is the identity.
If $\Aut_i(Y)/Z$ and $\Aut(L_0(Y))$
have the same dimension and number of connected components,
then any automorphism of $L_0(Y)$ is given as conjugation by an
element of $\Aut_i(Y)$.
\end{lem}

\begin{pf}
Let $\Ad : \Aut_i(Y)\to \Aut(L(Y))$ be given by $\Ad(g)(x) = gxg^{-1}$. 
From \cite{chevalley}, Proposition 12, Chapter III, No 9, it follows that 
the Lie algebra of the kernel of $\Ad$ is equal to the centre of $L(Y)$, 
which is spanned by $I_{n+1}$. Hence $\ker \Ad=Z$. Therefore the 
induced homomorphism $\Ad : \Aut_i(Y)/Z \to \Aut(L_0(Y))$ is injective. Lemma 
\ref{lem:surjective} implies that $\Ad$ is surjective, hence the
second statement follows. 
\end{pf}

For the three varieties that we consider we will check the hypothesis
of Lemma \ref{lem:lift} in separate subsections. 

\begin{rem} \rm
We remark that the hypothesis of Lemma \ref{lem:lift} does not hold for 
$Y=\p^2$. In order to construct the module isomorphism, one sometimes
has to correct the Lie algebra isomorphism by an outer automorphism
of $L_0(Y)$. See \cite{ghps} for details.
\end{rem}

For each of the three choices of $Y$, namely $\p^1$, $\PxP$, and $\mathbb{Y}$,
we still have to do three things:
\begin{enumerate}
\item solve the Lie algebra identification problem;
\item solve the Lie module identification problem;
\item check the hypothesis of Lemma~\ref{lem:lift}.
\end{enumerate}

Note that 2., the Lie module identification problem, is a just linear
problem, because the sought matrix $M$ as a generic solution
of the system of equations $M(xv)=x(Mv)$ for all $x\in L_0(Y)$ and $v\in F^{n+1}$.
A generic solution will be non-singular.
However, if the Lie algebra 
is split semisimple, then we can use the theory of weight vectors
to find a module isomorphism. This is more efficient then solving
the system of linear equations above.

\subsection{Identifying $\p^1$} \label{ss:p1}

Solutions for the identification problem for $\p^1$ are well-known.
Nevertheless, we want
to describe how solve it by Lie algebras because of two reasons:
first, it is the simplest possible example where the method works,
and second, we have to solve the corresponding Lie algebra identification
problem anyway at another place.

Using the anticanonical embedding,
we can reduce to the embedded identification problem of the parabola
with equation $y_0y_2-y_1^2=0$ in $\p^2$. The embedded twists 
of the parabola are exactly the nonsingular conics, and such a twist
is projectively isomorphic to the parabola iff it has a point defined
over $F$. Hence we see that our problem is
equivalent to deciding whether a given ternary quadratic form
is isotropic, i.e. has a nontrivial solution over $F$; 
constructing an explicit isomorphism is possible when we have 
such an explicit solution. 
We will see that the Lie algebra method reduces to the same problem.

\medskip

{\noindent{\em Identification of the Lie algebra.}}
We have that $\Aut_i(Y)/Z$ (where $Z$ consists of the scalar matrices)
is isomorphic to $\PGL_2(F)$. Therefore
$L_0(Y)\cong\ssl_2$. The twists of $\ssl_2$ are the
semisimple Lie algebras of dimension~3, because dimension and semisimplicity
do not change under field extension, and over algebraically closed fields
$\ssl_2$ is the only semisimple Lie algebra of dimension~3.
For checking semisimplicity, we can use Cartan's criterion
saying that this is equivalent to the Killing form being
non-degenerate. Finally, here is a proposition that allows
to identify a twist.

\begin{prop}
Let $L$ be a semisimple Lie algebra of dimension~3.
Then $L$ is isomorphic to $\ssl_2$ iff its Killing form
is isotropic.
\end{prop}

\begin{pf}
It is easy to check that the Killing form of $\ssl_2$ is isotropic,
hence ``only if'' is clear. 

Conversely, let $a\in L$ be a non-zero isotropic element. 
Note first, that for any nonzero $b$ in a twist of $\mf{sl}_2$ 
we have that the trace of $\ad(b)$ equals $0$ and also that 
the kernel of $\ad(b)$ is generated by $b$, for if it was
two-dimensional, then $[L,L]$ would be a nontrivial ideal in $L$.

Let $e_1,e_2,e_3$ be the eigenvalues of $\ad(a)$, so $e_1+e_2+e_3=0$.
Since $a$ is an isotropic element of the Killing form, we have also
$e_1^2+e_2^2+e_3^2=0$. One of eigenvalues is zero and so we get
that all $e_i$ vanish.
So $\ad(a)$ is nilpotent and
hence there exists an element $b$ such that $[a,b]=a$. Then $\ad(b)$
has an eigenvalue of $-1$, hence it is split semisimple and $b$
generates a split Cartan subalgebra $H$. 
When we have $H$, an isomorphism to $\ssl_2$ can
be constructed explicitly (see \cite{gra6}).
\end{pf}

Solving ternary quadratic forms can be done over number fields in Magma:
first, we check for local solvability at all primes dividing the Hessian.
If the form is everywhere
solvable, then there is a solution in $F$ by the Hasse principle
(see \cite{Omeara:71}). The construction of the solution can be reduced
to solving a norm equation of a quadratic extension of $F$.
If $F=\Q$, then we use faster algorithms for finding a rational
point on a plane conic.

\medskip

{\noindent{\em Identification of the Lie module.}}
We show that the $\ssl_2$-module given by the isomorphism
$\ssl_2\to L_0(Y)$ is irreducible, of highest weight $(2)$.
Let $V$ be the $2$-dimensional vector space over $F$ with basis $\{v_0,v_1\}$.
Let $W=\Sym^2(V)$ with the basis $\{v_0^2,2v_0v_1,v_1^2\}$.
Let $\varphi:V\to W$ be defined by $\varphi(v) = v^2$. 
We write the coordinates of an element of
$W$ with respect to the basis above. Then the image of the induced
map $\varphi: \p(V)\to\p(W)$ is exactly $Y$.

Let $\GL_2(F)$ act naturally on $V$, i.e.~the vector with the 
coordinates $(s,t)$ is mapped by $g=(g_{ij})_{i,j=0}^1\in\GL_2(F)$ 
to the vector with the coordinates $(g_{00}s+g_{01}t, g_{10}s+g_{11}t)$.
This leads to the action of $\GL_2(F)$ on $W$ by $g\cdot vv' = (gv)(gv')$,
for $v,v'\in V$. By writing the matrix of elements of $\GL_2(F)$ with 
respect to the basis above we get a representation $\rho : \GL_2(F)\to 
\GL_3(F)$. We have $g\cdot \varphi(v) = \varphi(g\cdot v)$, and hence 
$\varphi(V)$ is fixed under the action of $\GL_2(F)$ on $W$. We have 
further $Y = \varphi(\p(V))$, therefore $\rho(\GL_2(F))\subseteq
\Aut(\varphi(V))=\Aut_i(Y)$. The kernel of $\rho$ consists of two matrices, 
$\pm I_2$, the identity in $\GL_2(F)$. The conclusion is that the
$\GL_2(F)$-module given by $\rho$ is isomorphic to $\Sym^2(V)$. Hence
the same holds for the corresponding modules of the Lie algebras. 

Using highest weight vectors, we can construct a module isomorphism. 
This isomorphism is unique up to scalar multiplication, because the 
module is irreducible.

\medskip

{\noindent{\em Checking the hypothesis of Lemma~\ref{lem:lift}.}}
We have that $L_0(Y)$ is isomorphic to $\ssl_2(F)$ and hence
its centre is $0$. As in Lemma \ref{lem:lift}, let $Z$ be the subgroup
of $\Aut_i(Y)$ consisting of scalar multiples of the identity. Then 
both groups $\Aut_i(Y)/Z$ and $\Aut(\ssl_2)$ are isomorphic 
to $\PGL_2(F)$ (for $\Aut(\ssl_2)$ see \cite{jac}, Chapter IX, Theorem 5), 
hence both have dimension ~3 and are connected. 
{\hbox{}\hfill $\Box$\vskip0.2cm}

\subsection{Identifying $\PxP$} \label{ss:pxp}

Because $\PxP$ is a Del Pezzo surface of degree~8 (in its anticanonical
embedding), this identification problem is an instance of our
parametrization problem. 

\medskip

{\noindent{\em Identification of the Lie algebra.}}
Let $Z$ be the subgroup of $\Aut_i(Y)$ consisting of the scalar matrices.
Then $\Aut_i(Y)/Z \cong \PGL_2(F)\times \PGL_2(F)\ltimes \Z/2\Z$. Hence
$L_0(Y)\cong \slxsl$. Let $L$ be a given Lie algebra, then to
decide isomorphism with $L_0(Y)$ we do the following. 
We check whether it is semisimple, decompose into
its simple components, and solve the identification problem
for $\ssl_2$ (as in subsection~\ref{ss:p1}) for the two components.
If the number of components differs from~2, then $L$ is not
isomorphic to $L_0(Y)$. An algorithm for decomposing semisimple
Lie algebras can be found in \cite{gra6}.

\medskip

{\noindent{\em Identification of the Lie module.}}
The anticanonical embedding $i:\PxP\to\p^8$ is given by
\[ (s_0{:}s_1;t_0{:}t_1)\mapsto(s_0^2t_0^2{:}s_0^2t_0t_1{:}s_0^2t_1^2{:}
  s_0s_1t_0^2{:}s_0s_1t_0t_1{:}s_0s_1t_1^2{:}s_1^2t_0^2{:}s_1^2t_0t_1{:}
  s_1^2t_1^2) . \]
In this case the $\slxsl$-module given by the isomorphism
$\slxsl \to L_0(Y)$ is irreducible of highest weight $(2,2)$.
This can be shown as it was done for the $\p^1$ case (Section
\ref{ss:p1}). In this case we use the map $\varphi : 
V\times W \to \Sym^2(V)\otimes \Sym^2(W)$, where $V$, $W$ are
$2$-dimensional. Then the projectivization of $\varphi(V\times W)$
is equal to $Y$. Here the first direct summand of $\slxsl$ acts on
$\Sym^2(V)$ and the second summand on $\Sym^2(W)$. Hence the full
algebra $\slxsl$ acts on the tensor product. This means that the
module is irreducible and of highest weight $(2,2)$.

Hence we can decide module equivalence
by checking irreducibility and computing the highest weight.
In the affirmative case, we can again construct a module isomorphism
by using highest weight vectors. It is unique up to scalar multiplication,
as in the previous case.

\medskip

{\noindent{\em Checking the hypothesis of Lemma~\ref{lem:lift}.}}
Because $Y$ is anticanonically embedded, we have $\Aut_i(Y)/Z=\Aut(Y)$,
which has dimension~6 because its Lie algebra has dimension 6. The normal
subgroup $\Aut_0(Y)$ has the same dimension, but it is a proper subgroup
because the automorphism interchanging the two product factors $\p^1$
does not preserve classes. Hence $\Aut_i(Y)/Z$ has at least 2 components.
On the other hand, the group of automorphism of $\slxsl$ is
a semidirect product of the group of inner automorphism and
the finite group of ``diagram automorphisms'' (see  \cite{jac}, \S IX.4).
The group of inner automorphism is connected of dimension~6,
and the group of diagram automorphisms is $\Z/2\Z$ as the Dynkin
diagram consists of two nodes and no edges. Hence $\Aut(\slxsl)$
has dimension~6 and 2 connected components. The centre of $L_0(Y)$ is
$0$. So as seen in the proof of Lemma \ref{lem:lift}, the homomorphism 
$\Ad: \Aut_i(Y)/Z\to \Aut(\slxsl)$ is injective. Therefore $\Aut_i(Y)/Z$
has exactly two components. 
{\hbox{}\hfill $\Box$\vskip0.2cm}

\medskip

{\noindent{\em Timings.}}
We implemented the algorithm in Magma.
The examples were constructed as follows.
We took the canonical $\PxP$ in $\p^8$ given by 20 binomials.
Then we generated a $9\times 9$ matrix containing
random integer numbers with absolute values up to a given maximal
number (this is written in the first column of Table~\ref{tab:PxP}).
We used this matrix as the matrix of a linear transformation of 
projective space obtaining so a different system of implicit equations.
\begin{table}[!htb]
  \begin{center}
  \begin{tabular}{|r|r|r|r|r@{.}l|r@{.}l|r@{.}l|}
    \hline
    perturb & eqns max & LA size & prm size & \multicolumn{2}{|c|}{time} & 
    \multicolumn{2}{|l|}{LA time} & \multicolumn{2}{|l|}{conic time}\\
    \hline
    1       & 4        & 11      & 18       & 4&56   & 4&49   & 0&00 \\
    5       & 73       & 47      & 70       & 21&93  & 21&66  & 0&03 \\
    10      & 255      & 55      & 84       & 28&46  & 28&11  & 0&09 \\
    50      & 5026     & 84      & 130      & 48&75  & 48&15  & 0&22 \\
    100     & 25304    & 111     & 166      & 61&02  & 60&15  & 0&34 \\
    300     & 225440   & 134     & 200      & 75&86  & 73&00  & 2&14 \\
            & 208199   & 136     & 204      & 89&52  & 73&15  & 15&77 \\
    400     & 335499   & 143     & 213      & 77&99  & 76&31  & 0&93 \\
            & 418185   & 141     & 210      & 152&21 & 77&91  & 73&56 \\
            & 545728   & 140     & 208      & 482&69 & 74&50  & 407&53 \\
    500     & 720193   & 147     & 222      & 91&11  & 82&24  & 8&10 \\
            & 525179   & 145     & 216      & 80&95  & 78&91  & 1&29 \\
            & 546787   & 143     & 218      & 176&13 & 78&51  & 96&96 \\
    \hline
  \end{tabular}
  \end{center}
  \small
  \begin{tabular}{r@{ -- }p{11.3cm}}
    perturb & maximum entry allowed in perturbation matrix,\\
    eqns max & the maximal absolute value of the coefficients 
               in the implicit equations,\\
    LA size & the maximal length of the numerator/denominator 
              of the structure constants of the Lie algebra,\\
    prm size & the maximal length of the numerator/denominator 
               of the coefficients in the parametrization,\\
    time & the time (in sec) needed for parametrizing,\\
    LA time & the time (in sec) needed for finding the Lie 
              algebra (is a part of ``time'' in the previous column).\\
    conic time & the time (in sec) needed for finding rational points
                 on two conics constructed to identify two summands 
                 $\mf{sl}_2(\Q)$ (is a part of ``time''). \\
  \end{tabular}
  \caption{Parametrizing $\PxP$.}\label{tab:PxP}
\end{table}
For a ``small'' perturbation, almost the whole time is spent for 
finding the Lie algebra of the surface. As the coefficients of 
the linear transformation grow, finding a rational point on the conic 
starts to play the main role in the time complexity.

\subsection{Identifying $\mathbb{Y}$}

Because $\mathbb{Y}$, in its anticanonical embedding, 
is also a Del Pezzo surface of degree~8,
this identification problem is another instance of our
parametrization problem.

\medskip

{\noindent{\em Identification of the Lie algebra.}}

Every automorphism of $\mathbb{Y}$ leaves the exceptional line
invariant, so $\Aut(\mathbb{Y})$ is isomorphic to the subgroup
$\Aut(\p^2)=\PGL_3(F)$ fixing the point $(1{:}0{:}0)$. The whole
group leaves $\Pic(\mathbb{Y})$ invariant, so 
$\Aut_0(\mathbb{Y})=\Aut(\mathbb{Y})$. Its Lie algebra is isomorphic to
\[ L_0(\mathbb{Y})=\left\{ \begin{pmatrix} 2a & b_1 & b_2 \\ 
        0 & -a+c_1 & c_2 \\
        0 & c_3 & -a-c_1 \end{pmatrix} \mid a,b_1,b_2,c_1,c_2,c_3\in F 
        \right\} . \]
Here is a useful characterization of this Lie algebra.

\begin{prop}
Let $L$ be a Lie algebra. Then $L$ is isomorphic to $L_0(\mathbb{Y})$
iff it has a 2-dimensional ideal $I$ which is abelian as a subalgebra, 
and a 4-dimensional subalgebra $S$ isomorphic to $\gl_2$, such that
the adjoint action of $S$ on $I$ is faithful.
\end{prop}

\begin{pf}
``Only if'': for $L=L_0(\mathbb{Y})$, we take $I$ as the ideal defined by
$a=c_1=c_2=c_3=0$, and $S$ as the subalgebra defined by $b_1=b_2=0$.

``If'': the Lie algebra $\Der(I)$ is isomorphic to $\gl_2$. Because
any injective homomorphism from $\gl_2$ to itself is an automorphism,
the action of $S$ on $I$ is determined up to isomorphism. Therefore
$L$ is isomorphic to the semidirect sum $I\rtimes S$ with respect
to this action.
\end{pf}

To solve the identification problem for $L_0(\mathbb{Y})$ with
input $L$, we can proceed as follows.

\begin{enumerate}
\item Take $I$ as the nilradical of $L$.
        If this is not two-dimensional abelian, then
        $L$ is not isomorphic to $L_0(\mathbb{Y})$.
\item Take $S$ as the normalizer of a Levi subalgebra of $L$. 
        If $\dim(S)\ne 4$,
        then $L$ is not isomorphic to $L_0(\mathbb{Y})$.
\item Check if the adjoint action of $S$ on $I$ is faithful.
        If not, then $L$ is not isomorphic to $L_0(\mathbb{Y})$.
        If yes, one can construct an isomorphism using the
        construction of semidirect sums.
\end{enumerate}

For checking the correctness of the construction, it suffices to check
it for $L_0(\mathbb{Y})$; and this is a routine calculation.

\medskip

{\noindent{\em Identification of the Lie module.}}

Let $K$ be a Levi subalgebra of $L_0(\mathbb{Y})$ (for instance
the subalgebra defined by $a=b_1=b_2=0$). 
The given $L_0(\mathbb{Y})$-module $W=F^9$ is also
an $K$-module. We analyze this module by a similar method as the one
we used in Section \ref{ss:p1}.

Let $V$ be a $3$-dimensional vector space with basis $v_0,v_1,v_2$. 
Consider the symmetric power $\Sym^3(V)$ with the basis
$v_0^3$, $3v_0^2v_1$, $3v_0^2v_2$, 
$3v_0v_1^2$, $6v_0v_1v_2$, $3v_0v_2^2$, 
$v_1^3$, $3v_1^2v_2$, $3v_1v_2^2$, $v_2^3$.
Let $\varphi^\prime:V\to \Sym^3(V)$ be given by 
$\varphi^\prime(v) = v^3$.

Let $G=\GL_3(\Q)$ act naturally on $V$.
Let $\rho^\prime(g)$ be the matrix describing the action of $g\in G$ 
on $\Sym^3(V)$ with respect to the basis above.

Let $U$ be the subspace of $\Sym^3(V)$ spanned by $v_0^3$.
Let $\pi:\Sym^3(V)\to\Sym^3(V)/U$ be the projection discarding 
the coordinate at $v_0^3$, and set $\varphi = \pi\circ\varphi^\prime$.
Then $\mathbb{Y}$ is the projectivization of $\varphi(V)$.
Then $\Aut(\mathbb{Y})=\Stab_G(U)$. 

\begin{lem}\label{le:decomp}
  As a $K$-module $W$ decomposes as a direct sum
  $W = W_2 \oplus W_3 \oplus W_4$,
  where $W_i$ is an $i$-dimensional irreducible $K$-module.
  The elements of the nilradical $I$ carry $W_4$ to
  $W_3$ and $W_3$ to $W_2$.
\end{lem}

\begin{pf}
  When restricting to the Levi subalgebra $K$, 
  the module $\Sym^3(V)$ (see the discussion before the Lemma) 
  becomes an $\mf{sl}_2$-module and as such
  decomposes as a sum of four irreducible modules:
  $W_1 = U$,
  $W_2$ is the module spanned by $3v_0^2v_1, 3v_0^2v_2$ 
  and isomorphic to the natural $\mf{sl}_2$-module,
  $W_3$ is spanned by $3v_0v_1^2, 6v_0v_1v_2, 3v_0v_2^2$
  and isomorphic to $\Sym^2(F^2)$, and lastly
  $W_4$ is spanned by $v_1^3, 3v_1^2v_2, 3v_1v_2^2, v_2^3$
  and isomorphic to $\Sym^3(F^2)$.
  Then $W$ as $\mf{sl}_2$-module decomposes into the sum 
  $W_2 \oplus W_3 \oplus W_4$.

  To prove the last assertion of the Lemma,
  let $b\in I$, $b=b_1e_{12} + b_2e_{13}$, where
  $e_{ij}$ is the matrix with $1$ on the position $(i,j)$
  and $0$ elsewhere.
  So if $w\in W_4$ is a basis vector, $w=v_1^iv_2^{3-i}$, then 
  $b\cdot w \in \left< v_0v_1^{i-1}v_2^{3-i}, v_0v_1^iv_2^{2-i}\right>\subset W_3$.
  Similarly for $w\in W_3$ one gets $b\cdot w\in W_2$.
\end{pf}

Let $f : 
W\to W$ be an isomorphism of $L_0(\mathbb{Y})$-modules. Then $f$ restricted
to $W_i$ is multiplication by a scalar $\lambda_i$. Let $b=e_{12}\in I$,
and $w_4=v_1^3\in W_4$. Then $b\cdot v_1=v_0$, hence $b\cdot w_4 = 
3v_0v_1^2\in W_3$. Hence $f(b\cdot w_4)=\lambda_3b\cdot w_4$. On the
other hand, $f(b\cdot w_4) = b\cdot f(w_4)=\lambda_4 b\cdot w_4$. We infer
that $\lambda_4=\lambda_3$. In the same way we find that $\lambda_3=\lambda_2$,
so that $f$ is multiplication by a scalar. 

Now to identify the module we first decompose it into a direct sum
of irreducible $K$-modules. We note that this is straightforward using 
weight vectors. Then we find an isomorphism to $W$ by acting with 
elements of $I$, as in the discussion above. Again we have that
such an isomorphism is unique up to scalar multiplication.

\medskip

{\noindent{\em Checking the hypothesis of Lemma~\ref{lem:lift}.}}
The group $\Aut_0(\mathbb{Y})$ is connected and has dimension~6.
It suffices to prove that the automorphism group of $L_0(\mathbb{Y})$
is also connected and 6-dimensional. 

Any automorphism of $L_0(\mathbb{Y})$ is also an automorphism 
of the 3-di\-men\-sional radical $J$ 
and an automorphism of the two-dimensional nilradical $I$. 
The group of automorphisms of $I$ is $\GL_2(F)$, which is connected 
of dimension~4.
The subgroup of automorphisms of $J$ fixing $I$ pointwise is isomorphic
to $F^2$: an element in $x\in J-I$ can be mapped to any element in 
$y$ iff their adjoint action in $I$ is the same, and this is true
iff $x-y\in I$. Hence the group of automorphisms of $J$ is of dimension~6
and connected.
Finally, we show that any automorphism $\phi$ of $J$ can be extended 
in a unique way to an automorphism $\psi$ of $L(\mathbb{Y})$. 
Let $x\in J-I$ arbitrary. There is a unique Levi subalgebra $R$
that normalizes $x$. The automorphism $\psi$, if exists, has to
send $R$ to the unique subalgebra $R'$ that normalizes $\phi(x)$.
For any $y\in R$, there is a unique element $y'\in R'$ such that
$[y,z]=[y',z]$ for all $z\in I$. We set $\psi(y):=y'$, and this
determines the isomorphism $\psi$ uniquely. - It follows that
$\Aut(L_0(\mathbb{Y}))$ is isomorphic to $\Aut(J)$, hence it is 
also 6-dimensional and connected.
{\hbox{}\hfill $\Box$\vskip0.2cm}

\begin{rem} \label{rem} \rm
The algorithm for identifying $\mathbb{Y}$ does not require factorization
of polynomials or solving nonlinear equations; field arithmetic and
solving linear systems are sufficient. Hence the result---in particular
whether $L$ is isomorphic to $L_0(\mathbb{Y})$ or not---does not change
when we extend the field $F$. We rediscovered the well-known
fact that there are no proper twists of $\mathbb{Y}$ (see \cite{manin}).
\end{rem}

\medskip

{\noindent{\em Timings.}}
We tried our algorithm on examples which we constructed from the 
canonical surface (given by the binomial ideal with 20 generators) 
by a linear transformation of the projective space. 
The randomly generated matrix of the transformation has integral
entries with the given maximal absolute value
(the first column in Table~\ref{tab:blowup}).
We see that almost the whole time is spent for finding the Lie algebra
of the surface.
\begin{table}[htb]
  \begin{center}
  \begin{tabular}{|r|r|r|r|r@{.}l|r@{.}l|}
    \hline
    perturb & eqns max & LA size & prm size & 
    \multicolumn{2}{|c|}{time} & \multicolumn{2}{|l|}{LA time} \\
    \hline
    1       & 4        & 10      & 46       & 4&43   & 4&23 \\
    5       & 85       & 47      & 211      & 21&25  & 20&76 \\
    10      & 280      & 59      & 266      & 28&21  & 27&58 \\
    20      & 912      & 72      & 327      & 35&94  & 35&16 \\
    50      & 6372     & 93      & 424      & 51&66  & 50&43 \\
    100     & 26625    & 103     & 475      & 58&08  & 56&84 \\
    200     & 98407    & 127     & 584      & 69&29  & 67&69 \\
    500     & 599186   & 145     & 666      & 82&81  & 80&89 \\
    1000    & 1926906  & 159     & 724      & 91&26  & 89&11 \\
    2000    & 7973589  & 179     & 819      & 101&04 & 98&50 \\
    5000    & 60259495 & 207     & 957      & 118&99 & 115&94 \\
    10000   & 246171712& 219     & 1008     & 129&49 & 126&24 \\
    \hline
  \end{tabular}
  \end{center}
  \small
  Description of the columns: as in Table~\ref{tab:PxP}.
  \caption{Parametrizing $\mathbb{Y}$.}\label{tab:blowup}
\end{table}

\section{Parametrizing Twists of $\PxP$}\label{sec:simple}

The only Del Pezzo surfaces over algebraically closed fields are $\PxP$
and $\mathbb{Y}$. Hence any Del Pezzo surface over $F$ is a twist of
one of these two. There are no proper twists of $\mathbb{Y}$ by
remark~\ref{rem}, but
we still have to deal with proper twists of $\PxP$. (We will see that
some of them do have a parametrization.)

Here is a theorem that says that many twists do not have a parametrization.

\begin{thm} \label{thm:no}
Assume that $X\cong C_1\times C_2$, 
where $C_1$ and $C_2$ are twists of $\p^1$.
Then $X$ has a parametrization only if $C_1\cong\p^1$ and $C_2\cong\p^1$.
\end{thm}

\begin{pf}
Assume that $X$ has a parametrization. Then it has in particular
an $F$-rational point $p\in X(F)$.
The two projections give $F$-rational points $\pi_1(p)\in C_1(F)$ and 
$\pi_2(p)\in C_2(F)$. Because a twist of $\p^1$ with an $F$-rational
point is already isomorphic to $\p^1$, it follows that
$C_1\cong\p^1$ and $C_2\cong\p^1$.
\end{pf}

By Theorem~\ref{thm:no}, we can restrict our attention to varieties
that are not products. But how is this reflected in the Lie algebra?
Here is the answer to this question.

\begin{thm} \label{thm:2l}
A twist of $\PxP$ is a product of two twists of $\p^1$ iff its Lie
algebra is a direct sum of two twists of $\mf{sl}_2$.
\end{thm}

\begin{pf}
``Only if'': if $X\cong C_1\times C_2$, then $\Aut_0(X)$ is the direct
product of the two normal subgroups $\Aut(C_1)$ and $\Aut(C_2)$. It follows
that $L_0(X)=L_0(C_1)\oplus L_0(C_2)$.

``If'': assume that $X$ is not a product.
Let $E$ be a Galois extension of $F$ with the property that
$X_E\cong \PxP$. Then $\Pic(X_E)\cong\Z^2$, and the divisor classes
$(1,0)$ and $(0,1)$ define the two projections to $\p^1$. We claim that
the Galois group $G$ interchanges these two classes. Indeed, the action
of $G$ is $\Z$-linear, preserves the intersection product and the canonical
class $(-2,-2)$, and this shows that $(1,0)$ can only be mapped to itself
or to $(0,1)$. If $(1,0)$ was fixed, then the $G$-orbit sum of 
some divisor $D\in\Div(X_E)$ such that $[D]=(1,0)$ would be in 
$(|G|,0) = |G|(1,0)$,
and since it is in $\Div(X)$, 
it would then define a projection
to a twist of $\p^1$, contradicting our assumption that $X$ is not a product.

Since $G$ interchanges the two classes defining the two projections,
it also interchanges the two normal subgroups of $\Aut_0(X_E)$ of dimension~3,
and hence it also interchanges the two ideals of $L_0(X,E)$.
It follows that these ideals are not defined over $F$, hence $L_0(X,F)$
is simple.
\end{pf}

For any $a\in F^\ast-(F^\ast)^2$, we will now construct a twist $S_a$ 
of $\PxP$, called {\em sphere}, which is not a product, in the
simplest possible way.
More precisely,
let $E$ be the quadratic field extension $F[\alpha]/(\alpha^2-a)$.
Then $(S_a)_E\cong\PxP$.

The construction works as follows. We start with the anticanonical embedding
of $\p^1_E\times\p^1_E\subset\p^8_E$. We label coordinates and
unit vectors in $E^9$ by ordered pairs of integers in $\{0,1,2\}$. 
The surface $\PxP$ is embedded by mapping $((s:1),(t:1))$ to the
point with coordinates $x_{ij}=s^it^j$ with respect to the basis $e_{ij}$ 
for $i,j=0,1,2$.

Let $\sigma$ be the generator of the Galois group $G:={\cal G}(E,F)$. 
Then $\sigma$ induces an $F$-linear map $\Sigma:E^9\to E^9$ defined by 
$ce_{ij}\mapsto\sigma(c)e_{ji}$. Obviously $\Sigma$ preserves $\p^1_E\times\p^1_E$.
Similarly as in~\cite{delaunay}, the involution $\Sigma$ defines
an $F$-structure on $\p^1_E\times\p^1_E$. We set $S_a$ to be the $F$-variety
defined by this structure. The set of $F$-rational points on 
$\p^1_E\times\p^1_E$ is equal to the set of $E$-rational points 
fixed under $\Sigma$.

The variety $S_a$ is not a product because two factors
in $\PxP$ are interchanged by $\Sigma$, hence none of the two projection
morphisms is defined over $F$.

Let $V$ be the $F$-linear subspace of $E^9$ of fixed vectors.
By Galois descent, $\dim(V)=9$; we give the explicit basis
\begin{multline*}
 B := \{e_{00},e_{11},e_{22},e_{01}+e_{10},e_{12}+e_{21},e_{02}+e_{20}, \\
 \alpha^{-1}(e_{10}-e_{01}),
 \alpha^{-1}(e_{21}-e_{12}),\alpha^{-1}(e_{20}-e_{02}) \} .
\end{multline*}
We can give a parametrization of $S_a$ in the coordinates with respect
to the basis $B$ in the parameters $u:=\frac{1}{2}(s+t)$ and $v:=\frac{\alpha}{2}(s-t)$,
namely
\[ (1 : P : P^2 : u : Pu : 2u^2-P : v : vP : 2uv) , \]
where $P = u^2-a^{-1}v^2 = st$.

\subsection{Yet Another Identification Problem}

In this subsection, we give an algorithm for solving the embedded
identification problem for $S_a$. We denote its Lie algebra $L_0(S_a,F)$ 
by $\mf{s}_a$. We will show that
it is the $F$-linear space of elements in $\ssl_2(E)\oplus\ssl_2(E)$ that
are fixed under the semilinear automorphism that exchanges two fixed
Chevalley bases of the two summands and takes the coefficients to their
conjugates.

Of course, the algorithm first needs to find an $a\in F^\ast$ 
such that the given surface $X$ is isomorphic to $S_a$, if exists.

The {\em centroid} $\Gamma(L)$ of a Lie algebra $L$ is the centralizer
of $\ad L$ in $\gl(L)$. It is easy to check that the centroid of
$\mf{s}_a$ is isomorphic to $E:=F[\alpha]/(\alpha^2-a)$, the field
extension defined by $a$.

\begin{prop} \label{pr:cent}
Let $X$ be a twist of $\PxP$ which is not a product. Then the centroid $E$
of $L_0(X,F)$ is a quadratic field extension of $F$, and $X_E$ is a product.
\end{prop}

\begin{pf}
By Theorem~\ref{thm:2l}, we can assume that $L:=L_0(X,F)$ is simple.
By~\cite{jac}, Theorem~10.1, the centroid of a simple Lie algebra is a field.
Because $\Gamma(\slxsl)$ has dimension~2, and the dimension of the
centroid does not change when we extend the field, it follows
that $E=\Gamma(L)$ is a quadratic field extension. Because
$\Gamma(L\otimes_FE)=\Gamma(L)\otimes_FE=E\otimes_FE$ is not a field,
it follows that $L\otimes_FE$ is not simple.
\end{pf}

Of course, proposition~\ref{pr:cent} solves the subtask of finding $a$.
We just have to compute the centroid. Once we have $a$, there is of course
still no guarantee that $X$ is isomorphic to $S_a$; the following
proposition decides this.

\begin{prop} \label{pr:is}
Let $X$ be a twist of $\PxP$ which is not a product. 
Let $E:=F[\alpha]/(\alpha^2-a)$ be the centroid of $L_0(X,F)$.
Then the following are equivalent.

a) The varieties $X$ and $S_a$ are isomorphic.

b) The Lie algebras $L_0(X,F)$ and $\mf{s}_a$ are isomorphic.

c) The varieties $X_E$ and $\PxP$ are isomorphic over $E$.

d) The Lie algebras $L_0(X,E)$ and $\slxsl$ are isomorphic over $E$.
\end{prop}

\begin{pf}
(a)$\implies$(c): by construction, $(S_a)_E=\PxP$.

(c)$\implies$(d) (and also (a)$\implies$(b)) are obvious.

%
(d)$\implies$(b): in the following by $\sigma$-semilinear homomorphism $f$ 
of Lie algebras $L$, $L^\prime$ we mean an $F$-linear Lie algebra homomorphism
such that $f(cv) = \sigma(c)f(v)$ for every $c\in E$ and $v\in L$.
The Galois automorphism $\sigma$ induces a
$\sigma$-semilinear Lie algebra homomorphism $\sigma_L$ on
$L_0(X,E)=L_0(X,F)\otimes_FE$ which fixes $L_0(X,F)$.
By assumption, $L_0(X,E)$ is isomorphic to $\slxsl$, hence it is a sum of two
ideals $L_1$ and $L_2$, each isomorphic to $\ssl_2(E)$. The automorphism
$\sigma_L$ interchanges $L_1$ and $L_2$, because otherwise both would be
fixed under the Galois action and $L_0(X,F)$ would not be simple.
Let us fix a Chevalley basis in $\ssl_2(E)$ and let $\sigma_{\ssl_2}$
be the $\sigma$-semilinear automorphism of $\ssl_2$ fixing this basis.
Let $\psi:L_1\to\ssl_2(E)$ be a Lie algebra isomorphism. We define
the $E$-linear Lie algebra isomorphism $\varphi:L_1\oplus L_2\to\slxsl$ 
componentwise by sending $x_1\oplus x_2$ to 
$\psi(x_1)\oplus(\sigma_{\ssl_2}\circ\psi\circ\sigma_L)(x_2)$.
Let $\Sigma:\slxsl(E)\to\slxsl(E)$ be the semilinear automorphism
that interchanges the Chevalley bases of the summands,
i.e.~$\Sigma(x_1\oplus x_2) = \sigma_{\ssl_2}(x_2)\oplus\sigma_{\ssl_2}(x_1)$
Then the two semilinear Lie algebra homomorphisms
$\Sigma\circ\varphi$ and $\varphi\circ\sigma_L$ from $L_0(X,E)$ to
$\ssl_2(E)\oplus\ssl_2(E)$ coincide. It follows that the restriction
of $\varphi$ to $L_0(X,F)$ (as the subset of $L_0(X,E)$ which is fixed under
$\sigma$) is a Lie algebra isomorphism to the subset of
$\ssl_2(E)\oplus\ssl_2(E)$ the image of which is fixed under $\Sigma$,
and this is $\mf{s}_a$.

(b)$\implies$(a): the Lie algebra $\mf{s}_a$ acts on $F^9$ in two ways,
namely as the Lie algebra of $S_a$, and via the Lie algebra isomorphism
to $L_0(X,F)$ which we assume to exist. Over $E$, these two Lie modules
are both isomorphic to the unique irreducible module with highest weight
$(2,2)$ (see subsection~\ref{ss:pxp}). In particular, they are
isomorphic to each other. The matrix of a mo\-du\-le isomorphism 
describes also a Lie algebra isomorphism by conjugation.
Therefore it is a solution to a linear system and hence
defined over $F$.
Then by Section~\ref{sec:id} the claim follows.
\end{pf}

Here is the identification algorithm for $S_a$ applied to a given twist $X$
of $\PxP$ such that the centroid of $L_0(X,F)$ is $E:=F[\alpha]/(\alpha^2-a)$.

\begin{enumerate}
\item Decompose $L_0(X,E)$ into $L_1\oplus L_2$, using the algorithm
        described in \cite{gra6}.
\item Construct a Lie algebra isomorphism $\psi:L_1\to\ssl_2(E)$,
        using the algorithm described in subsection~\ref{ss:p1}.
        If the two Lie algebras are not isomorphic, then $X$
        is not isomorphic to $S_a$.
\item Construct a Lie algebra isomorphism $\varphi:L_0(X,F)\to\mf{s}_a$
        by restricting the $E$-isomorphism
        $\psi\oplus(\sigma\psi\sigma):L_1\oplus L_2\to\ssl_2(E)\oplus\ssl_2(E)$
        to $L_0(X,F)$.
\item Construct a Lie module isomorphism $M$ between the two 
        $\mf{s}_a$-modules given by the action on $S_a$ and by the 
        Lie algebra isomorphism $\varphi$. Return~$M$.
\end{enumerate}

The correctness of the algorithm follows from the following statements.

\begin{itemize}
\item Theorem~\ref{thm:2l} and Proposition~\ref{pr:cent} together
        imply that $L_0(X,E)$ decomposes into two ideals.
\item Proposition~\ref{pr:is} implies that $L_1\cong\ssl_2(E)$
        is necessary for $X$ being isomorphic to $S_a$.
\item The proof of Proposition~\ref{pr:is}, implication (d)$\implies$(b),
        shows that the construction in step~3 is indeed a Lie
        algebra isomorphism (hence $L_1\cong\ssl_2(E)$ is also sufficient
        for $X$ being isomorphic to $S_a$).
\item The proof of Proposition~\ref{pr:is}, implication (b)$\implies$(a),
        shows that the module isomorphism exists, is unique up to scalar
        multiplication, and takes $S_a$ into $X$.
\end{itemize}

\medskip

{\noindent{\em Timings.}}
For testing the algorithm we constructed examples as follows.
We have chosen $d\in\Z$ such that $d\not\in\Q^2$
(given in the first column of Table~\ref{tab:sphere}).
Then the sphere in $\p^3$ given by $z_0^2 - z_1^2 = z_2^2 - dz_3^2$
is isomorphic to $\PxP$ over $\Q(\sqrt{d})$ but not over $\Q$.
We anticanonically embedded the sphere over $\Q$ into $\p^8$ 
obtaining such a surface described by 14 binomials and 6 polynomials
with 4 terms. Afterwards we made a linear transformation similar 
to the two previous cases, just here the generated matrix is 
sparser, to obtain examples solvable in practice. 
Since we have to identify two $\mf{sl}_2$'s over $\Q(\sqrt{d})$,
we have to solve two relative norm equations.
This is very time consuming, therefore we were able to parametrize
only ``small'' examples.
\begin{table}[!htb]
  \begin{center}
  \begin{tabular}{|r|r|r|r|r|r@{.}l|r@{.}l|r@{.}l|}
    \hline
     & perturb & eqns & LA & prm & 
    \multicolumn{2}{|l|}{} & \multicolumn{2}{|c|}{LA} & 
    \multicolumn{2}{|c|}{normeq} \\
    discr & (sparse) & max & size & size & 
    \multicolumn{2}{|c|}{time} & \multicolumn{2}{|c|}{time} & 
    \multicolumn{2}{|c|}{time} \\
    \hline
    -1  & 1  & 3   & 3   & 9    & 2&460   & 0&670   & 1&010 \\
    3   & 1  & 5   & 3   & 23   & 3&620   & 1&030   & 1&700 \\
    8   & 1  & 15  & 5   & 1135 & 211&340 & 1&270   & 123&230 \\
    -1  & 2  & 10  & 5   & 92   & 41&690  & 1&250   & 38&670 \\
    \hline
  \end{tabular}
  \end{center}
  \small
  \begin{tabular}{r@{ -- }p{11.3cm}}
    discr & square of the primitive element used for the construction,\\
    normeq time & the time (in sec) needed for solving two relative 
                  norm equations (is a part of ``time''). \\
  \end{tabular}
  Description of the other columns: as in Table~\ref{tab:PxP}.
  \caption{Parametrizing the sphere.}\label{tab:sphere}
\end{table}

\subsection{Completeness of the Method}

Assume that $X$ is a twist of $\PxP$, which is not isomorphic to $\PxP$
and not isomorphic to $S_a$ for any $a\in F^\ast$. We distinguish two cases.

\begin{enumerate}
\item Assume that $X$ is a product. Then $X$ does not have a parametrization
        by Theorem~\ref{thm:no}.
\item Assume that $X$ is not a product. Let $E$ be the centroid of $L_0(X,F)$,
        which is a quadratic field extension by Proposition~\ref{pr:cent}.
        By Proposition~\ref{pr:is}, $X_E$ is not isomorphic to $\PxP$.
        On the other side, $X$ is a product by Proposition~\ref{pr:cent}.
        Then $X$ does not have a parametrization over $E$
        by Theorem~\ref{thm:no}. Consequently $X$ does not have a 
        parametrization over $F$.
\end{enumerate}

\bibliographystyle{plain}

\begin{thebibliography}{10}

\bibitem{magma}
Wieb Bosma, John Cannon, and Catherine Playoust.
\newblock The {M}agma algebra system. {I}. {T}he user language.
\newblock {\em J. Symbolic Comput.}, 24:235--265, 1997.
\newblock Computational algebra and number theory (London, 1993).

\bibitem{chevalley}
Claude Chevalley.
\newblock {\em Th\'eorie des groupes de {L}ie. {T}ome {II}. {G}roupes
  alg\'ebriques}.
\newblock Actualit\'es Sci. Ind. no. 1152. Hermann \& Cie., Paris, 1951.

\bibitem{ghps}
W.~A. de~Graaf, M.~Harrison, J.~P\'{\i}lnikov\'a, and J.~Schicho.
\newblock A {L}ie algebra method for rational parametrization of
  {S}everi-{B}rauer surfaces.
\newblock {\em Journal of Algebra}, 2005.

\bibitem{gra6}
Willem~A. de~Graaf.
\newblock {\em Lie algebras: theory and algorithms}, volume~56 of {\em
  North-Holland Mathematical Library}.
\newblock North-Holland Publishing Co., Amsterdam, 2000.

\bibitem{delaunay}
Claire Delaunay.
\newblock Real structures on smooth compact toric surfaces.
\newblock In {\em Topics in algebraic geometry and geometric modeling}, volume
  334 of {\em Contemp. Math.}, pages 267--290. Amer. Math. Soc., Providence,
  RI, 2003.

\bibitem{jac}
Nathan Jacobson.
\newblock {\em Lie algebras}.
\newblock Dover Publications Inc., New York, 1979.
\newblock Republication of the 1962 original.

\bibitem{manin}
Yu.~I. Manin.
\newblock {\em Cubic forms}, volume~4 of {\em North-Holland Mathematical
  Library}.
\newblock North-Holland Publishing Co., Amsterdam, second edition, 1986.
\newblock Algebra, geometry, arithmetic, Translated from the Russian by M.
  Hazewinkel.

\bibitem{Omeara:71}
O.~T. O'Meara.
\newblock {\em Introduction to quadratic forms}.
\newblock Springer-Verlag, New York, 1971.
\newblock Second printing, corrected, Die Grundlehren der mathematischen
  Wissenschaften, Band 117.

\bibitem{Schenck:04}
Hal Schenck.
\newblock Lattice polygons and {G}reen's theorem.
\newblock {\em Proc. Amer. Math. Soc.}, 132(12):3509--3512 (electronic), 2004.

\bibitem{josef1}
J.~Schicho.
\newblock Proper parametrization of surfaces with a rational pencil.
\newblock In {\em Proc. ISSAC 2000}, pages 292--299. ACM Press, 2000.

\bibitem{Shepperd-Barron:92}
N.~I. Shepperd-Barron.
\newblock The rationality of quintic {Del Pezzo} surfaces -- a short proof.
\newblock {\em Bull. London Math. Soc.}, 24:249--250, 1992.

\end{thebibliography}

\end{document}